\documentclass[12pt]{article}
\usepackage{amssymb}
\usepackage{amsmath}
\usepackage{amsthm}
\usepackage{cite}
\makeatletter
\def\LaTeX{\leavevmode L\raise.42ex
    \hbox{\kern-.3em\size{\sf@size}{0pt}\selectfont A}\kern-.15em\TeX}
\makeatother

\newcommand{\BibTeX}{{\rm B\kern-.05em{\sc
          i\kern-.025emb}\kern-.08em\TeX}}

\makeatletter
\def\@currentlabel{2.1}\label{e:dispaa}
\def\@currentlabel{2.21}\label{e:dispau}
\def\@currentlabel{2.22}\label{e:dispav}
\def\@currentlabel{2.23}\label{e:dispaw}
\def\@currentlabel{2.24}\label{e:dispax}
\def\theequation{\thesection.\@arabic\c@equation}
\makeatother

\newcounter{mnotecount}[section]

\newcommand{\rmnote}[1]{}

\renewcommand{\theequation}{\arabic{section}.\arabic{equation}}
\newtheorem{thm}{Theorem}[section]
\newtheorem{lem}{Lemma}[section]

\theoremstyle{definition}
\newtheorem{defn}{Definition}[section]
\newtheorem{rem}{Remark}[section]

\newcommand{\B}{\mathbb B}

\newcommand{\R}{\mathbb{R}}

\begin{document}
\title{Well-posedness of a fourth order evolution equation Modeling MEMS
\thanks {This research is supported by NSFC (11201119, 11471099)}}
\author{ Baishun Lai\\
\small {\it  Institute of Contemporary Mathematics, Henan University;}\\
\small {\it Kaifeng 475004, P.R.China}}
\date{}
\maketitle
\begin{abstract}
We consider a fourth order evolution equation involving a singular nonlinear term $\frac{\lambda}{(1-u)^{2}}$ in a bounded domain $\Omega\subset\R^{n}$.
This equation arises in the modeling of microelectromechanical systems. We first investigate the well-posedness of a fourth order parabolic equation which has been studied in \cite{Lau}, where the authors, by the semigroup argument, obtained the well-posedness of this equation for $n\leq2$. Instead of semigroup method, we use the Faedo-Galerkin technique to construct a unique solution of the fourth order parabolic equation for $n\leq7$, which improves and completes the result of \cite{Lau}. Besides, the well-posedness of the corresponding  fourth order hyperbolic equation is obtained by the similar argument for $n\leq7$.
\end{abstract}

\noindent
{\small {\bf Mathematics Subject Classification (2000):} 35J40, 35J75, 35J91.

\smallskip
\noindent
{\bf Key words:} Electrostatic MEMS, the fourth order evolution equation, well-posedness.}

\setcounter{equation}{0}
 \setcounter{equation}{0}
\section{Introduction}

Electrostatically actuated microelectromechanical systems (MEMS) are microscopic devices
which combine mechanical and electrostatic effects. MEMS devices have therefore become key components of
many commercial systems, including accelerometers for airbag deployment in automobiles,
ink jet printer heads, optical switches, chemical sensors, and so on (see, for
example, \cite{PB}).  A typical MEMS device is made of a rigid
conducting ground plate above which a clamped deformable plate (or membrane) coated with a thin
conducting film is suspended. An applied voltage difference between the two plates results in the deflection of the elastic plate, and a
consequent change in the MEMS capacitance, and thus transforms electrostatic energy into mechanical energy. The applied  voltage potential has an upper limit, beyond which the
electrostatic Coulomb force is not balanced by the elastic restoring force in the deformable plate, the two plates snap
together and the MEMS collapses. This phenomenon, called pull-in instability, was simultaneously observed
experimentally by Taylor \cite{T}, and Nathanson et al. \cite{N}. The critical displacement
and the critical voltage  potential associated with this instability are called pull-in displacement and pull-in voltage  potential,
respectively. Their accurate evaluation is crucial in the design of electrostatically actuated MEMS.

Mathematical models have been derived, see, for example, \cite{GP,LY,PB},  to describe the dynamics of the displacement $u=u(x,t)$ of the membrane $\Omega\subset \R^{n}$.  Let us sketch the derivation of this model  for the sake of completeness. Indeed, according to the Newton's second law and the narrow gap
asymptotic analysis, we see
$$
\gamma\frac{\partial^{2}u}{\partial t^{2}}=\mbox{electrostatic force}+\mbox{elastic force}+\mbox{damping force}
$$
where $\gamma$ is a constant denoted by the mass of membrane. Since we consider here the idealized situation where
the applied voltage and the permittivity of the membrane are constant (normalized
to one),  then
$$
\mbox{electrostatic force}=\lambda(\varepsilon^{2}|\nabla_{x}\psi(x,z,t)|^{2}+|\partial_{z}\psi(x,z,t)|^{2}),\  x\in \Omega,\  z>0,
$$
where $\lambda$ is proportional to the square of the applied voltage, $\psi$ is the the electrostatic potential and $\varepsilon$ denote the aspect ratio of the device. Under the small  aspect ratio condition
($\varepsilon\approx 0$), the $\psi$ is solved by
$$
\psi=\frac{(1-z)}{1-u},
$$
for details, see \cite{EP}.
Besides, we note that the damping force is linearly proportional to the velocity, that is
$$
\mbox{damping force}=-a\frac{\partial u}{\partial t},
$$
where $a$ is damping intensity,
and
$$
\mbox{elastic force}= \tau \Delta u-\beta \Delta^{2}u
$$
where $\tau$ is the tension constant in the stretching component of the energy, $\beta$ accounts for the bending energy. According to the
above discussion, the dimensionless dynamic deflection $u(x, t)$ of the membrane on a bounded domain $\Omega\subset \R^{n}$, under the small aspect assumption,  satisfies the following dynamic problem
\begin{align}\label{E1.1}
\left \{
\begin{array}{ll}
\gamma u_{tt}+a u_{t}+\beta \Delta^{2}u-\tau \Delta u=\frac{\lambda}{(1-u)^{2}},& x\in\Omega, t>0,\\
%0\leq u(x,t)< 1 ,& x\in\Omega, t>0,\\
u(x,0)=u^{0}(x),\ \ u_{t}(x,0)=u^{1}(x) & x\in\Omega,   \\
\mbox{boundary conditions},& x\in \partial\Omega, t>0.\\
 \end{array}
\right.
\end{align}
Observe that the right-hand side of equation features a singularity when $u=1$, which corresponds to the touchdown phenomenon already mentioned above.

The initial values $u^{0}(x), u^{1}(x) $ are  assumed to belong to some Sobolev space. Usually, one considers the following sets of boundary conditions
$$
u=\frac{\partial u}{\partial \nu}=0, \ \ \ x\in\partial\Omega,
$$
which we will refer to as Dirichlet boundary conditions, and
$$
u=\Delta u=0, \ \ \ x\in\partial\Omega,
$$
which we will refer to as Navier boundary conditions. The Dirichlet boundary condition is also called clamped boundary condition, which corresponds to the case where the capacitive actuator at the boundary is clamped, giving rise to zero vertical displacement and zero slope. Physically, the Navier boundary condition, usually referred to as the  pinned boundary condition, gives rise to a device which is ideally hinged along all its edges so that it is free to rotate and does not experience any torque or bending moment about its edges.

For the stationary case, \eqref{E1.1} has been studied extensively, see, for example, \cite{Co1,Gu4,Lai,LD}. For the non-stationary case, due to the lack of the maximum principle, little is known in the literature about the well-posedness of \eqref{E1.1} for $\beta>0$ so far. The author in \cite{Guo} established the local and global well-posedness of \eqref{E1.1} for pinned boundary conditions, $\gamma>0$ and the lower-dimensional case where $1\leq n\leq 3$.  Later, The authors in \cite{Lau} used the semigroup approach to obtain the existence of the  strong solutions of \eqref{E1.1}  for $n\leq 2, \gamma>0$. However, for the higher-dimensional case, the well-posedness of \eqref{E1.1} is open.

In the damping dominated limit $\gamma\ll 1$ when viscous forces dominate over inertial forces, \eqref{E1.1} reduces to the following forth order initial-boundary value parabolic problem

\begin{equation}\label{E1.2}
\left \{
\begin{array}{ll}
u_{t}+\beta \Delta^{2}u-\tau \Delta u=\frac{\lambda}{(1-u)^{2}},& x\in\Omega, t>0,\\
%0\leq u(x,t)< 1 ,& x\in\Omega, t>0,\\
u(x,0)=u^{0}(x), & x\in\Omega,   \\
\mbox{boundary conditions},& x\in \partial\Omega, t>0.\\
 \end{array}
\right.
\end{equation}
Here we let $a=1$ for simplicity.

In the present paper, we first investigate the local and global well-posedness of the parabolic  problem \eqref{E1.2}. When bending is neglected, that is, when $\beta=0$, this problem reduces to a second-order parabolic problem that has been studied extensively in the recent past, see, for example, \cite{EGG,GHW} and the references therein. Due to lack of the maximum principle, which plays an important role in studying the corresponding stationary problems,  only the references \cite{Lau, LL,LLS}, to the best our knowledge, give some partial results to this problem \eqref{E1.2} with $\beta>0$ so far.  To be more precise, the authors in \cite{Lau} use the semigroup argument to obtain the well-posedness of \eqref{E1.2} for any bounded domain $\Omega\subset\R^{n}$ and the lower-dimensional case $ n\leq 2$; the authors in \cite{LL,LLS}, by use of  numerical methods and asymptotic analysis,  considered the quenching phenomenon on a one-dimensional strip and the unit disc. In the present paper, we, instead of semigroup theory, will use the Faedo-Galerkin method to construct a solution of \eqref{E1.2} for $n\leq 7$, which improves and completes results of \cite{Lau}.

\begin{thm}\label{T1.1}
Let $\Omega\subset\R^{n}$ be an arbitrary bounded smooth domain for $n\leq 7$ and $\beta>0, \tau\geq 0, \lambda>0$.  Let
$u^{0}\in W^{4,2}(\Omega)\cap W_{0}^{2,2}(\Omega) $ be such that $\|u^{0}\|_{W^{4,2}(\Omega)\cap W_{0}^{2,2}(\Omega)}\leq \rho$ for some small $\rho\in (0,1)$. Then \eqref{E1.2} with Dirichlet boundary conditions  admits a unique solution $u(x,t)$ in
\begin{align*}
\mathcal{X}_{T}:=C^{0}([0,T];W^{4,2}(\Omega))\cap W^{1,2}(0,T; W_{0}^{2,2}(\Omega))\cap W^{1,\infty}(0,T;L^{2}(\Omega))
\end{align*}
with $\|u\|_{L^{\infty}(\Omega)}<1$, provided one of the following  conditions holds

(i) $\lambda\in \R^{+} $ and $T>0$ is sufficiently small;

(ii) $T=\infty$ and $\lambda\in \R^{+}$ is sufficiently small.

An identical result holds for the Navier problem but this time the solution belongs
to the space
$$
\mathcal{X}_{T}:=C^{0}([0,T];W^{4,2}(\Omega))\cap W^{1,2}(0,T; W_{0}^{1,2}(\Omega)\cap W^{2,2}(\Omega))\cap W^{1,\infty}(0,T;L^{2}(\Omega)).
$$

(iii) If $\lambda$ is sufficiently large, then $T_{m}<\infty$ for Dirichlet boundary conditions with $\Omega=\B_{1}$ (or Navier boundary conditions with any smooth domain $\Omega$). Here $\B_{1}$ is the unit ball and $T_{m}$ is the maximal existence time.
\end{thm}

\medskip
\begin{rem}
It is worth pointing out that the outcome of this Theorem complies with the physical viewpoint. More precisely, a ``pull-in'' instability occurs for high voltage values. Accordingly, for large values
of $\lambda$ solutions cease to exist globally, while solutions corresponding to small $\lambda$ values exist globally in time.
\end{rem}

\begin{rem}
For the third result on Dirichlet boundary conditions of the above Theorem,  a restriction with $\Omega=\B$ is needed. The essential reason of this is the lack of maximum principle in general domain.
\end{rem}

When inertial forces dominate over viscous forces in \eqref{E1.1}, i.e., $a\ll 1$, then \eqref{E1.1} reduces to the following hyperbolic problem (set $\gamma=1$ for simplicity)
\begin{equation}\label{E1.3}
\left \{
\begin{array}{ll}
u_{tt}+\beta \Delta^{2}u-\tau \Delta u=\frac{\lambda}{(1-u)^{2}},& x\in\Omega, t>0,\\
u(x,0)=u^{0}(x), u_{t}(x,0)=u^{1}(x) & x\in\Omega,   \\
u=\frac{\partial u}{\partial n}=0\ \  (\mbox{or}\ u=\Delta u=0) \ \ \ & x\in \partial\Omega,  t>0.\\
 \end{array}
\right.
\end{equation}
When $\beta=0$, this problem reduces to the second  hyperbolic problem which has been studied in \cite{K,LLZ}. For $\beta>0$
 this problem, to our knowledge, has not been investigated so far.
For this reason, we will give a result on its well-posedness though its argument is similar to the  parabolic case. To state our results precisely, we first introduce

\medskip
\begin{defn}
We call a function $u$ is weakly continuous from $[0,T]$ into the Banach space $Y$, if
$$
\forall\ v\in Y',  \mbox{the function}\ t\to <u(t), v> \mbox{is continuous},
$$
the set of all such functions will be denoted as $C_{w}([0,T]; Y)$.
\end{defn}

\medskip
\begin{thm}\label{T1.2}
Let $\Omega\subset\R^{n}$ be an arbitrary bounded smooth domain for $n\leq 7$ and $\beta>0, \tau\geq 0, \lambda>0$.  Let
$u^{0}\in W^{4,2}(\Omega)\cap W_{0}^{2,2}(\Omega), u^{1}\in W_{0}^{2,2}(\Omega)$  such that
 \begin{equation}\label{E1.4}
 \|u^{0}\|_{W^{4,2}(\Omega)}+\|u^{1}\|_{W^{2,2}(\Omega)}\leq \rho
 \end{equation}
 for some small $\rho\in (0,1)$.

(i) $\forall\ T>0, \exists\ \bar{\lambda}(T,\rho)>0$, if $0<\lambda<\bar{\lambda}(T,\rho)$, then \eqref{E1.3} with Dirichlet boundary conditions  admits a unique solution such that
\begin{align}
\begin{split}
u\in C_{w}(0,T; W^{4,2}(\Omega)), u'\in C_{w}(0,T; W^{2,2}_{0}(\Omega)),
u''\in L^{\infty}(0,T; L^{2}(\Omega)).
\end{split}
\end{align}

For the Navier problem: if the initial values $$u^{0}\in W^{4,2}(\Omega)\cap W_{0}^{1,2}(\Omega), u^{1}\in W^{2,2}(\Omega)\cap W_{0}^{1,2}(\Omega)$$ satisfies \eqref{E1.4}, then an identical result holds
 but this time the solution satisfies
\begin{align}
\begin{split}
u\in C_{w}(0,T; W^{4,2}(\Omega)), u'\in \  C_{w}(0,T; W^{2,2}(\Omega)\cap W^{1,2}_{0}(\Omega)),
u''\in L^{\infty}(0,T; L^{2}(\Omega)).
\end{split}
\end{align}

(ii) If $\lambda$ is sufficiently large, then the maximal existence time $T_{m}<\infty$ for Dirichlet boundary conditions with $\Omega=\B_{1}$ (or Navier boundary conditions with any smooth domain $\Omega$).
\end{thm}

Let us conclude this section with organization of the present paper as follows.

- in section 2 we recall some preliminary tools;

- in section 3 we will consider the well-posedness of the parabolic problem \eqref{E1.2}. To this end, we first study the well-posedness of the corresponding linear  parabolic problem which is of independent interested;

- in section 4 we will study the well-posedness of the hyperbolic problem \eqref{E1.3} by the same argument as that of section 3.

\setcounter{equation}{0}
\setcounter{equation}{0}
\section{Preliminaries}

Throughout the paper, we always suppose $\Omega\subset\R^{n}$ is a smooth bounded domain.
We denote by $\|\cdot\|_{p}$ the $L^{p}(\Omega)$ norm for $1\leq p\leq\infty$ and by $\|\cdot\|_{W^{s,p}}$ the $W^{s,p}(\Omega)$ norm.  Define
$$
(u,v)_{2}:=\int_{\Omega}uvdx \ \ \ \ \mbox{for all}\ \ u,v \in L^{2}(\Omega).
$$

On the space $W_{0}^{2,2}(\Omega)$, the bilinear form
\begin{equation}\label{E2.1}
(u,v)\mapsto (u,v)_{W_{0}^{2,2}}:=\beta(\Delta u,\Delta v)_{2}+\tau(\nabla u,\nabla v)_{2}\ \ \ \mbox{for all}\ \ u,v\in W_{0}^{2,2}(\Omega)
\end{equation}
define a scalar product over $W_{0}^{2,2}(\Omega)$ which induces a norm equal to $\|\cdot\|_{W^{2,2}(\Omega)}$. The space $W^{2,2}(\Omega)\cap W_{0}^{1,2}(\Omega)$ becomes a Hilbert space when endowed with the scalar product \eqref{E2.1},  please see \cite{GG} for details. Without loss of generality, we let $\mathcal{V}'$ denote the dual space of $\mathcal{V}$.

\medskip
\begin{lem}\label{L2.1}
(i) Each eigenvalue of $\mathcal{L}$ is real.

(ii) Furthermore, if we repeat each eigenvalue according to its (finite) multiplicity, all  the eigenvalues is given by
$$
\Sigma=\{\lambda_{k}\}_{k=1}^{\infty},
$$
where
$$
0<\lambda_{1}\leq \lambda_{2}\leq \lambda_{3}\leq\cdot\cdot\cdot,
$$
and
$$
\lambda_{k}\to\infty\ \ \ \mbox{as}\ \ k\to\infty.
$$

(iii) Finally, there exists an orthonormal basis $\{\omega_{k}\}_{k=1}^{\infty}$ of $L^{2}(\Omega)$, where $\omega_{k}\in W_{0}^{2,2}(\Omega)$ $( \mbox{or}\ W^{2,2}(\Omega)\cap W_{0}^{1,2}(\Omega))$ is an eigenfunction
corresponding to $\lambda_{k}$:
\begin{equation*}
\left \{
\begin{array}{ll}
\mathcal{L}\omega_{k}=\lambda_{k}\omega_{k}\ \ \ & \mbox{in}\ \ \Omega,\\
\omega_{k}=\frac{\partial\omega_{k}}{\partial n}=0\ \  (\mbox{or}\ \omega_{k}=\Delta\omega_{k}=0) \ \ \ & \mbox{on}\ \ \partial\Omega,
\end{array}
\right.
\end{equation*}
for $k=1,2,...$  Here $\mathcal{L}:=\beta\Delta^{2}-\tau \Delta$.
\end{lem}
\begin{rem}
By the regularity theory of the elliptic operator, $\omega_{k}\in C^{\infty}(\Omega)$ (and $\omega_{k}\in C^{\infty}(\bar{\Omega})$ if $\partial\Omega$ is smooth), for $k=1,2,...$.

\end{rem}
{\bf Proof.}\ \  By the Lax-Milgram theorem, $L^2$ theory of the elliptic operator and compact embedding theorem, we have
$$
\mathcal{S}:=\mathcal{L}^{-1}
$$
is bounded, linear, compact operator mapping $L^{2}(\Omega)$ into itself. Integrating by parts leads to
\begin{align*}
\int_{\Omega}g(\mathcal{S}f)dx=\int_{\Omega}\mathcal{L}(\mathcal{S}g)(\mathcal{S}f)dx=\int_{\Omega}\mathcal{L}(\mathcal{S}f)(\mathcal{S}g)dx=\int_{\Omega}f(\mathcal{S}g)dx,\ \ \ \forall f, g\in L^{2}(\Omega),
\end{align*}
which means that the operator $\mathcal{S}$ is self-adjoint. Therefore, by Hilbert-Schmidt's theorem, there exists a standard orthogonal basis $\{\omega_{k}\}_{k=1}^{\infty}$ such that
$\mathcal{S}\omega_{k}=\eta_{k}\omega_{k}, \eta_{k}\to0\ \ \ \mbox{as}\ \ k\to\infty$.  Notice also
$$
(\mathcal{S}f,f)_{2}=(u,f)_{2}=\beta(\Delta u,\Delta u)_{2}+\tau(\nabla u,\nabla u)_{2}\geq0 \ \  (f\in L^{2}(\Omega))
$$
and
\begin{equation*}
\left \{
\begin{array}{ll}
\mathcal{L}u=0\ \ \ & \mbox{in}\ \ \Omega,\\
u=\frac{\partial u}{\partial n}=0\ \  (u=\Delta u=0) \ \ \ & \mbox{on}\ \ \partial\Omega,
\end{array}
\right.
\end{equation*}
admits only a trivial solution, it is certainly $\eta_{k}>0$ and hence
\begin{equation*}
\left \{
\begin{array}{ll}
\mathcal{L}\omega_{k}=\frac{1}{\eta_{k}}\omega_{k}\ \ \ & \mbox{in}\ \ \Omega,\\
\omega_{k}=\frac{\partial \omega_{k}}{\partial n}=0\ \  (\omega_{k}=\Delta \omega_{k}=0) \ \ \ & \mbox{on}\ \ \partial\Omega.
\end{array}
\right.
\end{equation*}
The lemma follows. \qed

\medskip
The following Lemma is about the interpolation between $L^{2}(0,T;W^{m+4,2}(\Omega))$ and $W^{1,2}(0,T;W^{m,2}(\Omega))$.
\begin{lem}\label{L2.2}
Assume that $\Omega$ is open, bounded, and $\partial\Omega$ is smooth. Take $m$ to be a nonnegative integer. Suppose
$$
u\in L^{2}(0,T;W^{m+4,2}(\Omega)),\ \ \mbox{with}\ \ u'\in L^{2}(0,T;W^{m,2}(\Omega)),
$$
then
$$
u\in C([0,T];W^{m+2,2}(\Omega)),
$$
and
\begin{align}\label{E2.2}
\max_{0\leq t\leq T}\|u\|_{W^{m+2,2}(\Omega)}\leq C(\|u\|_{L^{2}(0,T;W^{m+4,2}(\Omega))}+\|u'\|_{L^{2}(0,T;W^{m,2}(\Omega))}).
\end{align}
\end{lem}
{\bf Proof.} The proof is standard, here we give a sketch of the proof. Suppose first that $m=0$, in which case
$$
u\in L^{2}(0,T; W^{4,2}(\Omega)), u'\in L^{2}(0,T; L^{2}(\Omega)).
$$
We select a bounded open set $\tilde{\Omega}\supset\supset \Omega$, and define a corresponding extension operator $E$ as follows:  $Eu=u$ a.e. in $\Omega$ and $Eu$ has support within $\tilde{\Omega}$.
We denote $Eu$ by $\bar{u}$ for simplicity. By the extension theorem of Sobolev space, we have
\begin{align}\label{E2.3}
\|\bar{u}\|_{L^{2}(0,T;W^{4,2}(\tilde{\Omega}))}\leq C \|u\|_{L^{2}(0,T;W^{4,2}(\Omega))};\ \
\|\bar{u}'\|_{L^{2}(0,T;L^{2}(\tilde{\Omega}))}\leq C \|u'\|_{L^{2}(0,T;L^{2}(\Omega))}.
\end{align}

We first claim that
$$
\bar{u}\in C([0,T];W^{2,2}(\tilde{\Omega})),
$$
from which we have
$$
u\in C([0,T];W^{2,2}(\Omega)).
$$
Indeed, we check as in the proof of Theorem 2 in section 5.3.2 of \cite{Ev} that there exist functions $\bar{u}^{\epsilon}(t)\in C^{\infty}(0,T;  W^{4,2}(\tilde{\Omega}))$ such that
as $\epsilon\to0$,
$$
\bar{u}^{\epsilon}(t)\to \bar{u}\ \ \ \mbox{in}\ \ L^{2}(0,T; W^{4,2}(\tilde{\Omega}))\cap W^{1,2}(0,T; L^{2}(\tilde{\Omega})).
$$
Now for $\epsilon,\delta>0$, we see that
\begin{align*}
\frac{d}{dt}\|\Delta(\bar{u}^{\epsilon}(t)-\bar{u}^{\delta}(t))\|^{2}_{2}&=2(\Delta(\bar{u}^{\epsilon}(t)-\bar{u}^{\delta}(t))',\Delta(\bar{u}^{\epsilon}(t)-\bar{u}^{\delta}(t)))\\&
=2((\bar{u}^{\epsilon}(t)-\bar{u}^{\delta}(t))',\Delta^{2}(\bar{u}^{\epsilon}(t)-\bar{u}^{\delta}(t)))
\end{align*}
Thus
\begin{align}\label{E2.4}
\begin{split}
\|\bar{u}^{\epsilon}(t)-\bar{u}^{\delta}(t)\|_{W^{2,2}(\tilde{\Omega}))}^{2}&=\|\bar{u}^{\epsilon}(s)-\bar{u}^{\delta}(s)\|_{W^{2,2}(\tilde{\Omega}))}^{2}\\&
+2\int_{s}^{t}((\bar{u}^{\epsilon}(\tau)-\bar{u}^{\delta}(\tau))',\Delta^{2}(\bar{u}^{\epsilon}(\tau)-\bar{u}^{\delta}(\tau)))d\tau
\end{split}
\end{align}
for all $0\leq s, t\leq T$. Fix any point $s\in (0,T)$ for which
$$
\bar{u}^{\epsilon}(s)\to \bar{u}(s)\ \ \  \mbox{in}\ \  W^{2,2}(\tilde{\Omega}).
$$
Then we have from \eqref{E2.4}
\begin{align*}
&\lim_{\epsilon,\delta\to0}\sup_{0\leq t\leq T}\|\bar{u}^{\epsilon}(t)-\bar{u}^{\delta}(t)\|_{W^{2,2}(\tilde{\Omega}))}^{2}\\
&\leq \lim_{\epsilon,\delta\to0}\int_{0}^{T}\Big(\|(\bar{u}^{\epsilon})'(\tau)-(\bar{u}^{\delta})'(\tau)\|_{2}^{2}+
\|\bar{u}^{\epsilon}(\tau)-\bar{u}^{\delta}(\tau)\|_{W^{4,2}(\tilde{\Omega})}^{2}\Big)d\tau\\&=0.
\end{align*}
i.e.,
$$
\bar{u}^{\epsilon}\to v\ \ \ \ \mbox{in}\ \ C([0,t]; W^{2,2}(\tilde{\Omega})).
$$
Besides, we also know that
$$
\bar{u}^{\epsilon}\to \bar{u}(t)\ \  \ \ \mbox{for}\  \mbox{a.e.}\ t
$$
and then $\bar{u}(t)=v(t)$ a.e. $t$. The claim follows.

Now we prove \eqref{E2.2}. Assume for the moment that $\bar{u}$ is smooth. We then compute
\begin{align*}
\frac{d}{dt}\int_{\tilde{\Omega}}|\Delta\bar{u}|^{2}dx=2\int_{\tilde{\Omega}}\Delta\bar{u}\Delta\bar{u}'dx=2\int_{\tilde{\Omega}}\Delta^{2}\bar{u}\bar{u}'dx
\leq(\|\bar{u}\|_{W^{4,2}(\tilde{\Omega})}^{2}+\|\bar{u}'\|_{2}^{2}).
\end{align*}
Thus
\begin{align}\label{E2.5}
\int_{\tilde{\Omega}}|\Delta\bar{u}(t)|^{2}dx\leq \int_{\tilde{\Omega}}|\Delta\bar{u}(s)|^{2}dx+
C(\|\bar{u}\|^{2}_{L^{2}(0,T;W^{4,2}(\tilde{\Omega}))}+\|\bar{u}'\|_{2}^{2})
\end{align}
for all $0\leq s, t\leq T$. We integrate \eqref{E2.5} with respect to $s$ and recall \eqref{E2.3} to obtain
\begin{align}\label{E2.6}
\max_{0\leq t\leq T}\|u\|_{W^{2,2}(\Omega)}\leq C(\|u\|_{L^{2}(0,T;W^{4,2}(\Omega))}+\|u'\|_{L^{2}(0,T;L^{2}(\Omega))}).
\end{align}
In the general case that $m>1$, we let $\alpha$ be a multiindex of order $|\alpha|\leq m$, and set $v:=D^{\alpha}u$. Then
$$
v\in L^{2}(0,T; W^{4,2}(\Omega)), v'\in L^{2}(0,T; L^{2}(\Omega)).
$$
We apply estimate \eqref{E2.6}, with $v$ replacing $u$, and sum over all indices $|\alpha|\leq m$ to obtain \eqref{E2.2}. We obtain the same estimate if $u$ is not smooth, upon approximating by
a smooth sequence $u^{\epsilon}$, as before.

\begin{lem}\label{L2.3}
Let $X$ and $Y$ be two Banach Spaces, such that $X\subset Y$ with a continuous injection. If a function $\phi$ belongs to $L^{\infty}(0,T;X)$ and is weakly continuous with values in $Y$, then $\phi$ is weakly continuous with values in $X$.
\end{lem}

For proof, please see \cite{T}, here we omit it.

\setcounter{equation}{0}
\setcounter{equation}{0}
\section{Well-posedness for the parabolic problem}

This section is devoted to the study of the parabolic problem \eqref{E1.2}. We first consider the well-posedness of the following  associated linear parabolic problem
\begin{equation}\label{E3.1}
\left \{
\begin{array}{ll}
u_{t}+\beta \Delta^{2}u-\tau \Delta u=f(x,t)\ \  &\mbox{in}\ \Omega\times (0,T),\\
u(x,0)=u^{0} \ \  &\mbox{in}\ \Omega,\\
u=\frac{\partial u}{\partial n}=0\ (\mbox{or}\ u=\Delta u=0) \ \  &\mbox{on}\ \Omega\times (0,T).
\end{array}\right.
\end{equation}

\begin{thm}\label{T3.1}
Let $0<\beta, 0<T\leq\infty$ and $f\in L^{2}(\Omega\times (0,T))$. The Dirichlet problem for the linear fourth order parabolic equation \eqref{E3.1}
with initial datum $u^{0}\in W_{0}^{2,2}(\Omega)$ admits a unique weak solution in the space
$$
C([0,T];W_{0}^{2,2}(\Omega))\cap L^{2}(0,T;W^{4,2}(\Omega))\cap W^{1,2}(0,T; L^{2}(\Omega)).
$$
The corresponding Navier problem with initial datum $u^{0}\in W^{2,2}(\Omega)\cap W_{0}^{1,2}(\Omega)$ admits a unique  weak solution in the space
$$
C([0,T];W^{2,2}(\Omega)\cap W_{0}^{1,2}(\Omega))\cap L^{2}(0,T;W^{4,2}(\Omega))\cap W^{1,2}(0,T; L^{2}(\Omega)).
$$
Furthermore, both cases admit the estimate
\begin{equation}\label{E3.2}
\max_{0\leq t\leq T}\| u(t)\|_{W^{2,2}(\Omega)}^{2}+\int_{0}^{T}\|u(t)\|_{W^{4,2}(\Omega)}^{2}+\int_{0}^{T}\|u_{t}\|_{2}^{2}\leq C(\|\Delta u^{0}\|_{2}^{2}
+\int_{0}^{T}\|f\|_{2}^{2})
\end{equation}
with the constant $C$ depending only on $\Omega, \beta, \tau$.
\end{thm}

\begin{defn}
We say a function
$$
u\in L^{2}(0,T; W_{0}^{2,2}(\Omega))\ \ (\mbox{or}\ L^{2}(0,T; W_{0}^{1,2}(\Omega)\cap W^{2,2}(\Omega)))$$
with
$$
 u'\in L^{2}(0,T; W^{-2,2}(\Omega))\ \  (\mbox{or}\ L^{2}(0,T;(W_{0}^{1,2}(\Omega)\cap W^{2,2}(\Omega))'))
$$
is a weak solution of the parabolic initial/ boundary-value problem \eqref{E3.1}  provided \medskip

(i) $$<u',v>+\beta(\Delta u, \Delta v)_{2}+\tau (D u, D v)_{2}=(f,v)_{2}$$
for each $v\in W_{0}^{2,2}(\Omega)\ (\mbox{or}\ W_{0}^{1,2}\cap W^{2,2}(\Omega))$ and a.e. time $0\leq t\leq T$, and \medskip

(ii) $u(0)=u^{0}$.
Here $<,>$ denotes the pairing between $W_{0}^{2,2}(\Omega)$ and $W_{0}^{-2,2}(\Omega)$ (or $W_{0}^{1,2}(\Omega)\cap W^{2,2}(\Omega)$ and $(W_{0}^{1,2}(\Omega)\cap W^{2,2}(\Omega))'$).
\end{defn}

\begin{rem}
In view of Lemma \ref{L2.2}, we see $u\in C([0,T];L^{2}(\Omega))$, and thus the equality (ii) makes sense.
\end{rem}

\bigskip
{\bf Proof Theorem \ref{T3.1}.}\ We will focus on Dirichlet boundary conditions, the proof for the Navier probelm follows with obvious modifications. Let
$u^{0}\in W_{0}^{2,2}(\Omega)$ and consider the following linear problem
\begin{equation}\label{E3.3}
\left \{
\begin{array}{ll}
u_{t}+\beta \Delta^{2}u-\tau \Delta u=f,& x\in\Omega, t>0,\\
u(x,0)=u^{0}(x), & x\in\Omega,   \\
u=\frac{\partial u}{\partial\nu}=0,& x\in \partial\Omega, t>0.\\
 \end{array}
\right.
\end{equation}

We intend to build a weak solution of \eqref{E3.3} by the so called {\em``Faedo-Galerkin''} method. More precisely,
let $\{\omega_{k}\}_{k=1}^{\infty}\subset W_{0}^{2,2}(\Omega)$ be an orthogonal complete system of eigenfunctions of $\beta \Delta^{2}-\tau \Delta $ under Dirichlet boundary conditions normalized by
$\|\omega_{k}\|_{2}=1$. By Lemma \ref{L2.1}
$$
\{\omega_{k}\}_{k=1}^{\infty}\ \mbox{is an orthonormal basis of }\ \ L^{2}(\Omega).
$$
Denote by $\{\lambda_{k}\}_{k=1}^{\infty}$ the unbounded sequence of corresponding eigenvalues. For any $k\geq1$ let
$$
u_{0}^{k}:=\sum_{i=1}^{k}(u^{0},\omega_{i})_{2}\omega_{i}
$$
so that $u_{0}^{k}\to u^{0}$ in $W_{0}^{2,2}(\Omega)$ as $k\to+\infty$. For each $k\geq 1$ we define an approximate solution $u_{k}:[0,T]\to W_{0}^{2,2}(\Omega)$ of \eqref{E3.3} as follows:
$$
u_{k}(t)=\sum_{i=1}^{k}g_{i}^{k}(t)\omega_{i}
$$
and
\begin{equation}\label{E3.4}
\left \{
\begin{array}{ll}
(u_{k}'(t),\omega_{j})_{2}+\beta (\Delta u_{k},\Delta\omega_{j} )_{2}+\tau (\nabla u_{k},\nabla \omega_{j})_{2}=(f(t),\omega_{j})_{2},\\
j=1,\cdot\cdot\cdot,k,\ \mbox{for a.e.}\ \ t\in (0,T),\\
u_{k}(x,0)=u^{k}_{0}(x),  x\in\Omega.
\end{array}
\right.
\end{equation}
So that for any $1\leq i\leq k$ the function $g_{i}^{k}(t)$ solves the Cauchy problem
\begin{equation}\label{E3.5}
\left \{
\begin{array}{ll}
(g_{i}^{k}(t))'+\sum_{i=1}^{k}\lambda_{i}g_{i}^{k}(t)=(f(t),\omega_{i})_{2}\\\\
g_{i}^{k}(0)=(u_{0}^{k},\omega_{i})_{2}.
\end{array}
\right.
\end{equation}
According to the standard existence theory for ordinary differential equations, the linear ordinary differential equation \eqref{E3.5} admits a unique solution $g_{i}^{k}$ such that $g_{i}^{k}\in W^{1,2}(0,T)$, and hence also \eqref{E3.4} admits
$u_{k}\in W^{1,2}(0,T;W_{0}^{2,2}(\Omega))$ as a unique solution.

We will obtain the a priori estimates independent of $k$ for the approximate solution $u_{k}$ and then pass to limit.

Step 1. {\em A priori estimates.} Indeed, we multiply the first equation of \eqref{E3.4} by
$g_{j}^{k}(t)$ and sum on $j$ from 1 up to $k$. We get
\begin{equation}\label{E3.6}
\frac{1}{2}\frac{d}{dt}\|u_{k}(t)\|_{2}^{2}+\beta \|\Delta u_{k}\|_{2}^{2}+\tau \|\nabla u_{k}\|_{2}^{2}=(f(t),u_{k})_{2}.
\end{equation}
Integrating over $(0,t)$  and using Cauchy's inequality with $\epsilon$, we are led to
$$
\|u_{k}(t)\|_{2}^{2}-\|u_{0}^{k}(t)\|_{2}^{2}+2\|u_{k}(t)\|_{L^{2}(0,T;W^{2,2}(\Omega))}^{2}
\leq\int_{0}^{T}\Big(4\|f(s)\|_{2}^{2}+\|u_{k}(s)\|_{W^{2,2}(\Omega)}^{2}\Big)ds.
$$
And therefore
\begin{align}\label{E3.7}
\begin{split}
\|u_{k}(t)\|_{L^{\infty}(0,T;L^{2}(\Omega))}^{2}+\|u_{k}(t)\|_{L^{2}(0,T;W^{2,2}(\Omega))}^{2}&\leq \|u_{0}^{k}(t)\|_{2}^{2}
+4\|f\|_{L^{2}(0,T;L^{2}(\Omega))}\\&
\leq 4(\|u^{0}\|_{2}^{2}+\|f\|_{L^{2}(0,T;L^{2}(\Omega))}).
\end{split}
\end{align}
 Next we multiply both sides of \eqref{E3.4} by $(g_{j}^{k}(t))'$ and sum on $j$ to obtain
$$
(u_{k}'(t),u_{k}'(t))_{2}+\beta (\Delta u_{k},\Delta u_{k}'(t) )_{2}+\tau (\nabla u_{k},\nabla u_{k}'(t))_{2}=(f(t),u_{k}'(t))_{2}
$$
Integrating over $(0,T)$  and using Cauchy's inequality with $\epsilon$, we see that
\begin{align}\label{E3.8}
\begin{split}
&\int_{0}^{T}\int_{\Omega}|u_{k}'(t)|^{2}dxdt+\frac{\beta}{2}\|\Delta u_{k}(\cdot,T)\|_{L^{2}(\Omega)}^{2}+\frac{\tau}{2}\|\nabla u_{k}(\cdot,T)\|_{L^{2}(\Omega)}^{2}\\&
\leq \frac{\beta}{2}\|\Delta u_{k}(\cdot,0)\|_{L^{2}(\Omega)}^{2}+\frac{\tau}{2}\|\nabla u_{k}(\cdot,0)\|_{L^{2}(\Omega)}^{2}+\int_{0}^{T}\int_{\Omega}|f|^{2}dxdt\\&
\leq \|u^{0}\|_{W^{2,2}(\Omega)}^{2}+\int_{0}^{T}\int_{\Omega}|f|^{2}dxdt.
\end{split}
\end{align}

Step 2. {\em Passage to limit.} From \eqref{E3.7} and \eqref{E3.8}, we may extract a subsequence, still denoted by $\{u_{k}\}$ such that
$$
u_{k}\rightharpoonup u\ \ \ \mbox{in}\ \  L^{2}(0,T;W^{2,2}(\Omega))\cap W^{1,2}(0,T; L^{2}(\Omega)).
$$
We expect that the limit function $u$ is a weak solution of \eqref{E3.3}. To this end, we introduce a function $h\in C^{1}([0,T]; C_{0}^{2}(\Omega))$ and
take a an approximate sequence of $h$
$$
h_{j}(x,t)=\sum_{m=1}^{j}\alpha_{j,m}(t)\omega_{m}(x)
$$
such that $\|h_{j}-h\|_{L^{2}(0,T; W_{0}^{2,2}(\Omega))}\to0$ as $j\to\infty$. Here $\{\alpha_{j,m}(t)\}_{j=1}^{m}$ are given smooth functions.  Now multiplying the first equation of (\ref{E3.4}) by
$\alpha_{j,m}$ and summing on $m$ from 1 up to $j$, we, by taking the limit for $k\to\infty$, see that
$$
\int_{0}^{T}(u'(t),h_{j})_{2}+\int_{0}^{T}\beta(\Delta u,\Delta h_{j} )_{2}+\int_{0}^{T}\tau(\nabla u, \nabla h_{j} )_{2}=\int_{0}^{T}(f(t),h_{j})_{2}.
$$
Letting $j\to\infty$ then we are led to
\begin{equation}\label{E3.9}
\int_{0}^{T}(u'(t),h)_{2}+\int_{0}^{T}\beta(\Delta u,\Delta h)_{2}+\int_{0}^{T}\tau(\nabla u, \nabla h)_{2}=\int_{0}^{T}(f(t),h)_{2}.
\end{equation}
Since $C^{1}([0,T]; C_{0}^{2}(\Omega))$ is dense in $L^{2}(0,T; W_{0}^{2,2}(\Omega))$,  we conclude equality \eqref{E3.9} is valid for any $h\in L^{2}(0,T; W_{0}^{2,2}(\Omega))$, which further implies
$$
(u'(t),h)_{2}+\beta(\Delta u,\Delta h)_{2}+\tau(\nabla u, \nabla h)_{2}=(f(t),h)_{2}
$$
for all $h\in W_{0}^{2,2}(\Omega)$ and a.e. $0\leq t\leq T$.

Now we claim $u(0)=u^{0}$. Indeed, from \eqref{E3.9} we deduce that
\begin{equation}\label{E3.10}
\int_{0}^{T}-(v',u)_{2}+\beta(\Delta u,\Delta v)_{2}+\tau(\nabla u, \nabla v)_{2}dt
=\int_{0}^{T} (f,v)_{2}dt+(u(0),v(0))_{2}
\end{equation}
for each $v\in C^{1}([0,T]; W_{0}^{2,2}(\Omega))$ with $v(T)=0$. Similar, from \eqref{E3.4} we also have
$$
\int_{0}^{T}-(v',u_{k})_{2}+\beta(\Delta u_{k},\Delta v)_{2}+\tau(\nabla u_{k}, \nabla v)_{2}dt
=\int_{0}^{T} (f,v)_{2}dt+(u_{k}(0),v(0))_{2}
$$
Let $k \to \infty$, we deduce that
\begin{equation}\label{E3.11}
\int_{0}^{T}-(v',u)_{2}+\beta(\Delta u,\Delta v)_{2}+\tau(\nabla u, \nabla v)_{2}dt
=\int_{0}^{T} (f,v)_{2}dt+(u_{0},v(0))_{2}
\end{equation}
Here we have used the fact $u_{k}(0)\to u^{0}$ in $L^{2}(\Omega)$. As $v(0)$ is arbitrary, comparing \eqref{E3.10} and \eqref{E3.11}, we conclude $u(0)=u^{0}$.
From this and \eqref{E3.9}, we conclude
$$u\in L^{2}(0,T;W^{2,2}(\Omega))\cap W^{1,2}(0,T; L^{2}(\Omega))$$
is a weak solution of \eqref{E3.3} which satisfies \eqref{E3.2}. Uniqueness follows from the contradiction argument:
if $v,w$ were two solutions of \eqref{E3.3} which share the same initial date, by subtracting the equations and \eqref{E3.2} we would get
\begin{equation*}
\max_{0\leq t\leq T}\|\Delta (v-w)\|_{2}^{2}+\int_{0}^{T}\|\Delta^{2}(v-w)\|_{2}^{2}+\int_{0}^{T}\|(v-w)_{t}\|_{2}^{2}\leq 0,
\end{equation*}
which immediately yields $v\equiv w$.

Step 3. {\em Ends of proof.}
Since
\begin{equation*}
\left \{
\begin{array}{ll}
\beta \Delta^{2}u=f-u_{t}+\tau \Delta u\in L^2(\Omega\times(0,T)),& x\in\Omega, t>0,\\
u=\frac{\partial u}{\partial\nu}=0,& x\in \partial\Omega, t>0,\\
 \end{array}
\right.
\end{equation*}
then we have $u\in L^{2}(0,T;W^{4,2}(\Omega))$ by the regularity theorem of elliptic operator. Taking advantaging of interpolation between
$L^{2}(0,t;W^{4,2}(\Omega))$ and $W^{1,2}(0,T; L^{2}(\Omega))$, we obtain $u\in C(0,T; W_{0}^{2,2}(\Omega)).$

\begin{thm}\label{T3.2}(Improved regularity).
If $u^{0}\in W^{4,2}(\Omega)\cap W^{2,2}_{0}(\Omega)\ (\mbox{or}\ W^{4,2}(\Omega)\cap W^{1,2}_{0}(\Omega)),  \ \ f'\in L^{2}(0,T; L^{2}(\Omega)),$ then
$$
u\in L^{\infty}(0,T; W^{4,2}(\Omega)), u'\in L^{\infty}(0,T; L^{2}(\Omega))\cap L^{2}(0,T; W_{0}^{2,2}(\Omega)),
$$
with the estimate
\begin{align}\label{E3.12}
\mbox{ess}\sup_{0\leq t\leq T}(\|u'(t)\|_{2}^{2}+\|u\|_{W^{4,2}(\Omega)}^{2})+\int_{0}^{T}\|u'(t)\|_{W^{2,2}(\Omega)}^{2}dt\leq C(\|f\|_{W^{1,2}(0,T;L^{2}(\Omega))}^{2}+\|u^{0}\|_{W^{4,2}(\Omega)}).
\end{align}
Here the constant $C$ depends only on $\Omega, \beta, \tau$.
\end{thm}

{\bf Proof.}\ Fix $k\geq 1$ and differentiate equation (\ref{E3.4}) with respect to $t$, we find
\begin{equation}\label{E3.13}
(\tilde{u}_{k}',\omega_{j})_{2}+\beta(\Delta\tilde{u}_{k},\Delta \omega_{j})_{2}+\tau(\nabla\tilde{u}_{k},\nabla\omega_{j})_{2}=(f', \omega_{j})_{2},\ \  (j=1,\cdot\cdot\cdot,k)
\end{equation}
where $\tilde{u}_{k}:=u_{k}'$. Multiply (\ref{E3.13}) by $\frac{d}{dt}g_{k}^{j}(t)$ and sum $j=1,\cdot\cdot\cdot,k$, we see
$$
(\tilde{u}_{k}',\tilde{u}_{k})_{2}+\beta(\Delta\tilde{u}_{k},\Delta\tilde{u}_{k})_{2}+\tau(\nabla\tilde{u}_{k},\nabla\tilde{u}_{k})_{2}=(f',\tilde{u}_{k})_{2}.
$$
Integrating over $(0,T)$  and using Cauchy's inequality with $\epsilon$, we deduce
\begin{align}\label{E3.14}
\begin{split}
\sup_{0\leq t\leq T}\|u_{k}'(t)\|_{2}^{2}&+2\beta\int_{0}^{T}\|\Delta u_{k}'(t)\|_{2}^{2}dt
+2\tau\int_{0}^{T}\|\nabla u_{k}'(t)\|_{2}^{2}dt\\&\leq C(\|u_{k}'(0)\|_{2}^{2}+\|f'\|^{2}_{L^{2}(0,T;L^{2}(\Omega))})\\&
\leq C(\|f\|^{2}_{W^{1,2}(0,T;L^{2}(\Omega))}+\|u_{k}(0)\|_{W^{4,2}(\Omega)}^{2}).
\end{split}
\end{align}
Here, we employed the first equation of (\ref{E3.4}) in the last inequality.

Remember that we have taken $\{\omega_{k}\}$ to be the complete collection of (smooth) eigenfunctions for $\beta\Delta^{2}-\tau\Delta$ on $W_{0}^{2,2}(\Omega)$. In particular
$(\beta\Delta^{2}-\tau\Delta) u_{k}=0$ on $\partial\Omega$. Thus
$$
\|u_{k}(0)\|_{W^{4,2}(\Omega)}^{2}\leq C\|(\beta\Delta^{2}-\tau\Delta) u_{k}(0)\|_{L^{2}}^{2}=C(u_{k}(0), (\beta\Delta^{2}-\tau\Delta)^{2}u_{k}(0))_{2}.
$$
Since $(\beta\Delta^{2}-\tau\Delta)^{2}u_{k}(0)\in \mbox{span}\{\omega_{j}\}_{j=1}^{k}$ and $(u_{k}(0), \omega_{j})_{2}=(u^{0},\omega_{j})_{2}$ for $j=1,\cdot\cdot\cdot,k$, we have
\begin{align*}
\|u_{k}(0)\|_{W^{4,2}(\Omega)}^{2}&\leq C(u^{0},(\beta\Delta^{2}-\tau\Delta)^{2}u_{k}(0))_{2}=C((\beta\Delta^{2}-\tau\Delta)u^{0},(\beta\Delta^{2}-\tau\Delta)u_{k}(0))_{2}\\&
\leq \frac{1}{2}\|u_{k}(0)\|_{W^{4,2}(\Omega)}^{2}+C\|u^{0}\|_{W^{4,2}(\Omega)}^{2}.
\end{align*}
Therefore, combining with (\ref{E3.14}), we have
\begin{align}\label{NE3.15}
\sup_{0\leq t\leq T}\|u_{k}'(t)\|_{2}^{2}+\int_{0}^{T}\|u_{k}'(t)\|_{W^{2,2}(\Omega)}^{2}dt\leq C(\|f\|_{W^{1,2}(0,T;L^{2}(\Omega))}^{2}+\|u^{0}\|_{W^{4,2}(\Omega)}^{2})
\end{align}

Now $$((\beta\Delta^{2}-\tau\Delta)u_{k},\omega_{j})_{2}=(f-u_{k}',\omega_{j})_{2}\ \ (j=1,\cdot\cdot\cdot,k).$$ And multiplying this identity by $\lambda_{j}g_{k}^{j}(t)$ and
summing $j=1,\cdot\cdot\cdot,k$, we deduce for $0\leq t\leq T$ that
$$
((\beta\Delta^{2}-\tau\Delta)u_{k},(\beta\Delta^{2}-\tau\Delta)u_{k})_{2}=(f-u_{k}',(\beta\Delta^{2}-\tau\Delta)u_{k})_{2},
$$
By the H\"{o}lder inequality, we see that
\begin{align*}
\|\Delta^{2}u_{k}\|_{2}^{2}&\leq(f-u_{k}'(t),(\beta\Delta^{2}-\tau\Delta)u_{k})_{2}+C\|u_{k}\|_{W^{2,2}(\Omega)}\\&
\leq C\|f\|_{2}^{2}+\frac{1}{2}\|\Delta^{2}u_{k}\|_{2}^{2}+C\|u_{k}'(t)\|_{2}^{2}+C\|u_{k}\|_{W^{2,2}(\Omega)}
\end{align*}
and then passing to limits as $k=k_{l}\to\infty$ and combining \eqref{E3.2} and \eqref{NE3.15}, we deduce
\begin{align*}
\sup_{0\leq t\leq T}(\|u'(t)\|_{2}^{2}+\|u\|_{W^{4,2}(\Omega)}^{2})+\int_{0}^{T}\|u'(t)\|_{W^{2,2}(\Omega)}^{2}dt\leq C(\|f\|_{W^{1,2}(0,T;L^{2}(\Omega))}^{2}+\|u^{0}\|^{2}_{W^{4,2}(\Omega)}).
\end{align*}

{\bf Proof of Theorem \ref{T1.1}.}\ \ We only consider Dirichlet boundary conditions, the proof for the Navier probelm is similar.
 Since we consider the case $1\leq N\leq 7$, then by the Sobolev embedding theorem, we deduce
\begin{align}\label{E3.15}
\|u\|_{L^{\infty}(\Omega)}\leq C(\Omega) \|u\|_{W^{4,2}(\Omega)}.
\end{align}
Now define
$$
\mathcal{X}_{T}:=C^{0}([0,T];W^{4,2}(\Omega))\cap W^{1,2}(0,T; W_{0}^{2,2}(\Omega))\cap W^{1,\infty}(0,T;L^{2}(\Omega))
$$
with norm
$$
\|v\|_{\mathcal{X}_{T}}^{2}:=\int_{0}^{T}(\|v_{t}\|_{W^{2,2}(\Omega)}^{2}+\|v\|_{W^{2,2}(\Omega)}^{2})dt+\max_{0\leq t\leq T}(\|v\|_{W^{4,2}(\Omega)}^{2}+\|v_{t}\|_{2}^{2}).
$$
And define
\begin{align*}
\bar{M}(R,T):=\{v\in  \mathcal{X}_{T}:  \|v\|_{\mathcal{X}_{T}}\leq R\}
\end{align*}
with $R$ satisfying $C(\Omega)R<1$. Here $C(\Omega)$ is defined in \eqref{E3.15}.
Let $0<r<R$, we also define the set
\begin{align*}
M(r,T):=\{v\in\mathcal{X}_{T}:  \|v\|_{\mathcal{X}_{T}}<r\}.
\end{align*}
 From \eqref{E3.15},  we have
\begin{equation} \label{E3.16}
u(t)\in \bar{M}(R,T)\Rightarrow \|u\|_{L^{\infty}(\Omega\times(0,T))}\leq C(\Omega)R<1.
\end{equation}

Now let $r\in (0, R)$ be fixed and
$$
u_{i}(t) \in\bar{M}(r,T),
$$ for $i=1,2$, then by the Theorem \ref{T3.1}, the initial-Dirichlet linear problem
\begin{equation}\label{E3.17}
\left \{
\begin{array}{ll}
v_{t}+ \beta \Delta^{2}v-\tau \Delta v=\frac{\lambda}{(1-u_{i})^{2}},& x\in\Omega, t>0,\\
v(x,0)=u^{0}(x), v_{t}(x,0)=u^{1}(x) & x\in\Omega,   \\
v=\frac{\partial v}{\partial\nu}=0,& x\in \partial\Omega, t>0,\\
 \end{array}
\right.
\end{equation}
has a unique solution $$v_{i}(t):=\mathcal{F}(u_{i})\in C([0,T];W_{0}^{2,2}(\Omega))\cap L^{2}(0,T;W^{4,2}(\Omega))\cap W^{1,2}(0,T; L^{2}(\Omega))$$ for $i=1,2$.

Now we claim $$v_{i}\in  C([0,T];W^{4,2}(\Omega)),$$
from which we have $v_{i}\in\mathcal{X}_{T}$.
Indeed, since  $\|u_{i}\|_{L^{\infty}(\Omega)}<1$ and $u_{i}\in W^{1,2}(0,T; L^{2}(\Omega))$, we have
$$\frac{\lambda}{(1-u_{i})^{2}}\in W^{1,2}(0,T; L^{2}(\Omega)).$$ And then by Theorem \ref{T3.2}, we see
$$
v_{i}\in W^{1,\infty}(0,T; L^{2}(\Omega))\cap W^{1,2}(0,T; W_{0}^{2,2}(\Omega)).
$$
From this, it is easy to see that
$$
\beta\Delta^{2}v_{i}=\frac{\lambda}{(1-u_{i})^{2}}+\tau\Delta v_{i}-\frac{dv_{i}(t)}{dt}\in W^{2,2}(\Omega), \mbox{a.e.}\  t.
$$
Then by the regularity theory of the elliptic operator,  we are led to
$$
v_{i}\in L^{2}(0,T; W^{6,2}(\Omega)).
$$
Combining $v_{i}\in W^{1,2}(0,T; W_{0}^{2,2}(\Omega))$ with Lemma \ref{L2.2}, we have
$$
v_{i}\in C([0,T]; W^{4,2}(\Omega)).
$$

Using the Theorem \ref{T3.2} again, we see that
\begin{align*}%\label{E2.20}
\begin{split}
\|v_{1}-v_{2}\|_{\mathcal{X}_{T}}&\leq\lambda C \|\frac{1}{(1-u_{1})^{2}}-\frac{1}{(1-u_{2})^{2}}\|_{W^{1,2}(0,T;L^{2}(\Omega))}\\
&=\lambda C\Big(\int_{0}^{T}\int_{\Omega}|\frac{1}{(1-u_{1})^{2}}-\frac{1}{(1-u_{2})^{2}}|^{2}dxdt\Big)^{\frac{1}{2}}\\&+2\lambda C\Big(\int_{0}^{T}\int_{\Omega}|\frac{u_{1}'}{(1-u_{1})^{3}}-\frac{u_{2}'}{(1-u_{2})^{3}}|^{2}dxdt\Big)^{\frac{1}{2}}\\
&=:I+II.
\end{split}
\end{align*}
For $I$, we have
\begin{align}\label{E3.18}
\begin{split}
I&= 2\lambda C \Big(\int_{0}^{T}\int_{\Omega}\frac{(u_{1}-u_{2})^{2}}{(1-(\theta u_{1}+(1-\theta)u_{2}))^{6}}dxdt\Big)^{\frac{1}{2}}
\\&\leq 2\lambda C k(r)\Big(\int_{0}^{T}\int_{\Omega}(u_{1}-u_{2})^{2}dxdt\Big)^{\frac{1}{2}}\\&
\leq 2\lambda C k(r)\|u_{1}-u_{2}\|_{\mathcal{X}_{T}}.
\end{split}
\end{align}
or
\begin{align}\label{E3.19}
I\leq 2\lambda C k(r)T^{\frac{1}{2}}\|u_{1}-u_{2}\|_{C([0,T];W^{4,2}(\Omega))}\leq 2\lambda C k(r)T^{\frac{1}{2}}\|u_{1}-u_{2}\|_{\mathcal{X}_{T}}.
\end{align}
For $II$, we have
\begin{align}\label{E3.20}
\begin{split}
II&\leq 2\lambda C \Big(\int_{0}^{T}\int_{\Omega}\frac{(u'_{1}-u'_{2})^{2}}{(1-u_{1})^{6}}dxdt\Big)^{\frac{1}{2}}\\&
+2\lambda C\Big(\int_{0}^{T}\int_{\Omega}|u'_{2}|^{2}\Big|\frac{1}{(1-u_{1})^{3}}-\frac{1}{(1-u_{2})^{3}}\Big|^{2}dxdt\Big)^{\frac{1}{2}}\\&
\leq k(r)2\lambda C\Big[\|u_{1}-u_{2}\|_{\mathcal{X}_{T}}+\|u_{1}-u_{2}\|_{L^{\infty}((0,t)\times\Omega)}\Big(\int_{0}^{T}\int_{\Omega}|u_{2}'|^{2}dxdt\Big)^{\frac{1}{2}}\Big]\\&
\leq2\lambda C(r+k(r))\|u_{1}-u_{2}\|_{\mathcal{X}_{T}}.
\end{split}
\end{align}
or
\begin{align}\label{E3.21}
\begin{split}
II&\leq 2\lambda CT^{\frac{1}{2}}k(r)\Big(\|u_{1}-u_{2}\|_{W^{1,\infty}(0,T;L^{2}(\Omega))}+\|u_{1}-u_{2}\|_{L^{\infty}((0,t)\times\Omega)}\|u_{2}\|_{W^{1,\infty}(0,T;L^{2}(\Omega))}\Big)\\&
\leq2\lambda C(r+k(r))T^{\frac{1}{2}}\|u_{1}-u_{2}\|_{\mathcal{X}_{T}}.
\end{split}
\end{align}
Here and in what follows $k(r)$ is a positive nondecreasing function for $r\in [0,R_{0}]$ and $C$  depends only on $\Omega$. From \eqref{E3.18} and \eqref{E3.20}, we have
\begin{align}\label{E3.22}
\|v_{1}-v_{2}\|_{\mathcal{X}_{T}}\leq 2\lambda C(r+k(r))\|u_{1}-u_{2}\|_{\mathcal{X}_{T}},
\end{align}
and from \eqref{E3.19} and \eqref{E3.21}
\begin{align}\label{E3.23}
\|v_{1}-v_{2}\|_{\mathcal{X}_{T}}\leq 2\lambda T^{\frac{1}{2}} C(r+k(r))\|u_{1}-u_{2}\|_{\mathcal{X}_{T}}.
\end{align}

Now consider the unique solution $w(t)$ to the linear problem
$$
w_{t}+\beta \Delta^{2}w-\tau \Delta w=0,\ \ \ x\in \Omega, t>0,
$$
with the same boundary and initial conditions as \eqref{E3.3}. By the Theorem \ref{T3.1}, we have
 $$w(t)\in C([0,T];W_{0}^{2,2}(\Omega))\cap L^{2}(0,T;W^{4,2}(\Omega))\cap W^{1,2}(0,T; L^{2}(\Omega))$$
such that
 \begin{align}
 \begin{split}
\|w\|_{\mathcal{X}_{T}}&\leq C\|u^{0}\|_{W^{4,2}(\Omega)}:=C\rho.
\end{split}
\end{align}

Define the ball
$$
\mathcal{B}_{\frac{r}{2}}=\{u\in \mathcal{X}_{T}: \|u-w\|_{\mathcal{X}_{T}}\leq \frac{r}{2}\}.
$$
Choosing $\rho$  small enough  such that
$$
C\rho+\frac{r}{2}<r,
$$
we then have, if $u(t)\in \mathcal{B}_{\frac{r}{2}}$,
$$
\|u(t)\|_{\mathcal{X}_{T}}\leq \|w(t)\|_{\mathcal{X}_{T}}+\frac{r}{2}<r,\ \ \forall\ 0\leq t\leq T,
$$
and then $\mathcal{B}_{\frac{r}{2}}\subset M(r,T)$. \medskip

Case 1. Global existence for small $\lambda$. Now using estimate \eqref{E3.22}, we find
\begin{align}\label{E3.25}
\|v_{i}-w\|_{\mathcal{X}_{T}}\leq 2C\lambda (k(r)+r)r
\end{align}
for $i=1,2$. Now choosing $\lambda$ so small that
$$
\lambda\leq \lambda(r):= \frac{1}{4C(k(r)+r)},
$$
we have from \eqref{E3.22} and \eqref{E3.25} that
\begin{align*}
\begin{split}
&\|v_{1}-v_{2}\|_{\mathcal{X}_{T}}\leq \frac{1}{2}\|u_{1}-u_{2}\|_{\mathcal{X}_{T}};\\
&\|v_{i}-w\|_{\mathcal{X}_{T}}\leq\frac{r}{2}.
\end{split}
\end{align*}
Hence the map
\begin{align}
\begin{split}
\mathcal{F}:& \mathcal{B}_{\frac{r}{2}}\to \mathcal{B}_{\frac{r}{2}}\\
& u_{i}\to v_{i}\ \ (i=1,2)
\end{split}
\end{align}
is a contraction map and it has a unique fixed point $u=\mathcal{F}(u)$ in $\mathcal{B}_{\frac{r}{2}}$ for $0<\lambda\leq\lambda(r)$ and arbitrary $T>0$, which is a global weak solution of \eqref{E1.2} with Dirichlet boundary conditions.

Case 2. Local existence in time. Similarly, using estimate \eqref{E3.2}, we are led to
\begin{align}\label{E3.27}
\|v_{i}-w\|_{\mathcal{X}_{T}}\leq C\lambda T^{\frac{1}{2}}(k(r)+r)\|u_{i}\|_{\mathcal{X}_{T}}
\leq \lambda T^{\frac{1}{2}} C(k(r)+r)r
\end{align}
for $i=1,2$. Let $T$ small enough such that
$$
0<T^{\frac{1}{2}}\leq \bar{T}(\lambda,\rho,r):=\frac{1}{2C\lambda(k(r)+r)},
$$
we then have from \eqref{E3.23} and \eqref{E3.27} that
\begin{align*}
\begin{split}
&\|v_{1}-v_{2}\|_{\mathcal{X}_{T}}\leq \frac{1}{2}\|u_{1}-u_{2}\|_{\mathcal{X}_{T}};\\
&\|v_{i}-w\|_{\mathcal{X}_{T}}\leq\frac{r}{2}.
\end{split}
\end{align*}
The existence of a unique solution to \eqref{E1.2} over $[0,T]$ for all $T\leq \bar{T}(\lambda,r)$ follows from the application of the Banach fixed point Theorem to the map.

Now we give the proof of (iii) of Theorem \ref{T1.1} as follows. To this end, we will use the eigenfunction method which comes from, for example, \cite{CL,Lac,Lau}. Indeed, from \cite{GG},  there exists a pair $(\lambda_{1},\phi_{1})$ such that
$0<\lambda_{1}, 0<\phi_{1}\in C^{4}(\bar{\B})\cap W_{0}^{2,2}(\B), \|\phi_{1}\|_{1}=1$ and
$$
\left \{
\begin{array}{ll}
\beta \Delta^{2}\phi_{1}-\tau \Delta \phi_{1}=\lambda_{1}\phi_{1},& x\in\B,\\
\phi_{1}=\frac{\partial \phi_{1}}{\partial\nu}=0,& x\in \partial\B.\\
 \end{array}
\right.
$$
Let $u(x,t)$ be the solution on $[0,T_{m})$ to \eqref{E1.2} and define for $t\in [0,T_{m})$
$$
M(t):=\int_{\B}\phi_{1}(x)u(x,t)dx\leq \int_{\B}\phi_{1}dx=1.
$$
Now we multiply \eqref{E1.2} by $\phi_{1}$, integrate over $\B_{1}$ and use the properties of $\phi_{1}$ and Jensen's inequality to obtain
\begin{align}\label{E3.28}
\begin{split}
\frac{dM}{dt}&=-\int_{\B}(\beta\Delta^{2}\phi_{1}-\tau\Delta\phi_{1})udx+\lambda\int_{\B}\frac{\phi_{1}}{(1-u)^{2}}dx\\&
\geq-\lambda_{1}\int_{\B}\phi_{1}udx+\frac{\lambda}{(1-\int_{\B}\phi_{1}udx)^{2}}\\&
=-\lambda_{1}M+\frac{\lambda}{(1-M)^{2}}:=g(M)
\end{split}
\end{align}
By a simple calculation, we have $g(M)>c_{0}>0$ if we choose $\lambda>\frac{4\lambda_{1}}{27}$. From \eqref{E3.28}, we immediately have
$$
1-M(0)\geq M(t)-M(0)\geq c_{0}t,
$$
consequently, $T_{m}\leq \frac{1-M(0)}{c_{0}}<\infty$.

\setcounter{equation}{0}
\setcounter{equation}{0}
\section{Well-posedness for the hyperbolic problem}

In this section, we will consider the well-posedness of the dynamic problem \eqref{E1.3}. As in Section 3, we first study the well-posedness of the corresponding the linear hyperbolic problem
\begin{equation}\label{E4.1}
\left \{
\begin{array}{ll}
u_{tt}+\beta \Delta^{2}u-\tau \Delta u=f(x,t) & x\in\Omega, t>0,\\
u(x,0)=u^{0}(x), u_{t}(x,0)=u^{1}(x) & x\in\Omega,   \\
u=\frac{\partial u}{\partial n}=0\ \  (\mbox{or}\ u=\Delta u=0) \ \ \ & x\in \partial\Omega,  t>0,\\
 \end{array}
\right.
\end{equation}
where $u^{0}(x), u^{1}(x) $ are  assumed to belong to some Sobolev space, $f(x,t)\in L^{2}(\Omega\times (0,T))$.

\begin{defn}
We say a function
$$
u\in L^{2}(0,T; W_{0}^{2,2}(\Omega))\ \ (\mbox{or}\ L^{2}(0,T; W_{0}^{1,2}(\Omega)\cap W^{2,2}(\Omega)))$$
with
$$
 u''\in L^{2}(0,T; W^{-2,2}(\Omega))\ \  (\mbox{or}\ L^{2}(0,T;(W_{0}^{1,2}(\Omega)\cap W^{2,2}(\Omega))'))
$$
is a weak solution of the hyperbolic initial/ boundary-value problem \eqref{E4.1}  provided \medskip

(i) $$<u'',v>+\beta(\Delta u, \Delta v)_{2}+\tau (\nabla u, \nabla v)_{2}=(f,v)_{2}$$
for each $v\in W_{0}^{2,2}(\Omega)\ (\mbox{or}\ W_{0}^{1,2}(\Omega)\cap W^{2,2}(\Omega))$ and a.e. time $0\leq t\leq T$, and \medskip

(ii) $u(0)=u^{0}, u'(0)=u^{1}$.
Here $<,>$ denotes the pairing between $W_{0}^{2,2}(\Omega)$ and $W^{-2,2}(\Omega)$ (or $W_{0}^{1,2}(\Omega)\cap W^{2,2}(\Omega)$ and $(W_{0}^{1,2}(\Omega)\cap W^{2,2}(\Omega))'$).
\end{defn}

\vskip0.2in

\begin{thm}\label{T4.1}
Let $0<T<\infty$ and $f\in L^{2}(\Omega\times (0,T))$. The Dirichlet problem for the linear fourth order hyperbolic equation \eqref{E4.1}
with initial datums $u^{0}\in W_{0}^{2,2}(\Omega), u^{1}\in L^{2}(\Omega)$ admits a unique weak solution such that
$$
u\in C_{w}([0,T];W_{0}^{2,2}(\Omega)), u'(t)\in C_{w}([0,T];L^{2}(\Omega)),  u''(t)\in L^{2}(0,T;W^{-2,2}(\Omega)).
$$
And the corresponding Navier problem with initial datums $u^{0}\in W^{2,2}(\Omega)\cap W_{0}^{1,2}(\Omega), u^{1}\in L^{2}(\Omega)$ admits a unique  weak solution such that
\begin{align*}
&u\in  C_{w}([0,T];W^{2,2}(\Omega)\cap W_{0}^{1,2}(\Omega)), u'(t)\in C_{w}([0,T];L^{2}(\Omega)), \\& u''(t)\in L^{2}(0,T;(W^{2,2}(\Omega)\cap W_{0}^{1,2}(\Omega))').
\end{align*}
Furthermore, both cases admit the estimate
\begin{equation}\label{E4.2}
\mbox{ess}\sup_{0\leq t\leq T}(\|u\|_{W^{2,2}(\Omega)}^{2}+\|u'(t)\|_{2}^{2})+\int_{0}^{T}\|u''\|^{2}_{W^{-2,2}(\Omega)}dt\leq C\Big(\|\Delta u^{0}\|_{2}^{2}
+\|u^{1}\|_{2}^{2}+\int_{0}^{T}\|f\|_{2}^{2}dt\Big).
\end{equation}
Here the constant $C$ depends only on $\Omega, T,\beta,\tau$.
\end{thm}
\vskip 0.2in
{\bf Proof of Theorem \ref{T4.1}.} As in Theorem \ref{T3.1}, we only consider  the Dirichlet boundary condition case, the proof for the Navier problem follows with the obvious modification.
Similar to Theorem \ref{T3.1}, we will once more employ {\em``Faedo-Galerkin''} method to construct our weak solutions. To this end, we, exactly as in the proof of Theorem \ref{T3.1}, define an approximate solution $u_{k}:[0,T]\to W_{0}^{2,2}(\Omega)$ of \eqref{E4.1} as follows:
$$
u_{k}(x,t)=\sum_{i=1}^{k}g_{i}^{k}(t)\omega_{i}(x),\ \ \ \ k\geq1,
$$
where $\omega_{i}(x)$ is defined as in Lemma \ref{L2.1} and the function $g_{i}^{k}(t)\ (1\leq i\leq k)$ solves the Cauchy problem
\begin{equation}\label{E4.3}
\left \{
\begin{array}{ll}
(g_{i}^{k}(t))''+\lambda_{i}g_{i}^{k}(t)=(f(t),\omega_{i})_{2},\\\\
g_{i}^{k}(0)=(u_{0}^{k},\omega_{i})_{2}, \frac{d}{dt}g_{i}^{k}(0)=(u_{1}^{k},\omega_{i})_{2},
\end{array}
\right.
\end{equation}
with
$$
u_{0}^{k}(x):=\sum_{i=1}^{k}(u^{0},\omega_{i})_{2}\omega_{i}(x); \ \  u_{1}^{k}(x):=\sum_{i=1}^{k}(u^{1},\omega_{i})_{2}\omega_{i}(x).
$$
According to the standard theory for ordinary differential equations, there exists a unique function $g_{i}^{k}(t)\in W^{2,2}(0,T)$ solving \eqref{E4.3} for $0\leq t\leq T$. \medskip

As in proof of Theorem \ref{T3.1}, we first study {\em a priori estimates} of the approximate solution $u_{k}$. Indeed,
\begin{equation}\label{E4.4}
(u_{k}''(t),\omega_{j})_{2}+\beta(\Delta u_{k},\Delta \omega_{j})_{2}+\tau(\nabla u_{k},\nabla \omega_{j})_{2}=(f,\omega_{j})_{2},
\end{equation}
multiply this equality by $\frac{d}{dt}g_{j}^{k}(t)$, sum $j=1,...,k$, we see
$$
(u_{k}''(t),u_{k}'(t))_{2}+\beta(\Delta u_{k},\Delta u_{k}')_{2}+\tau(\nabla u_{k},\nabla u_{k}')_{2}=(f,u_{k}')_{2}.
$$
From this, we immediately have
\begin{equation}\label{E4.5}
\frac{d}{dt}\Big(\|u_{k}'\|_{L^{2}}^{2}+\beta\|\Delta u_{k}\|_{L^{2}}^{2}+\tau\|\nabla u_{k}\|_{L^{2}}^{2}\Big)\leq C(\|u_{k}'\|_{L^{2}}^{2}
+\|f\|_{L^{2}}^{2}).
\end{equation}
Now write
$$
\eta(t):=\|u_{k}'\|_{L^{2}(\Omega)}^{2}+\beta\|\Delta u_{k}\|_{L^{2}(\Omega)}^{2}+\tau\|\nabla u_{k}\|_{L^{2}(\Omega)}^{2}.
$$
Then inequality \eqref{E4.5} reads
$$
\eta'(t)\leq C(\eta(t)+\|f\|_{L^{2}}^{2})
$$
for $0\leq t\leq T$. Thus Grownwall's inequality yields the estimate
$$
\eta(t)\leq e^{Ct}\Big(\eta(0)+\int_{0}^{t}\|f(s)\|_{L^{2}}^{2}ds\Big),
$$
where
\begin{align*}
\eta(0)&=\|u'_{k}(0)\|_{L^{2}}^{2}+\beta\|\Delta u_{k}(0)\|_{L^{2}}^{2}+\tau\|\nabla u_{k}(0)\|_{L^{2}}^{2}\\&
\leq C(\|u^{1}\|_{L^{2}(\Omega)}^{2}+\|u^{0}\|_{W^{2,2}(\Omega)}^{2}).
\end{align*}
Thus, we are led to
\begin{align}\label{E4.6}
\max_{0\leq t\leq T}\Big(\|u_{k}'\|_{L^{2}(\Omega)}^{2}+\|u_{k}\|_{W^{2,2}(\Omega)}^{2}\Big)\leq
C\Big(\|u^{1}\|_{L^{2}(\Omega)}^{2}+\|u^{0}\|_{W^{2,2}(\Omega)}^{2}+\int_{0}^{T}\|f(s)\|_{L^{2}}^{2}ds\Big).
\end{align}

Fix any $v\in W_{0}^{2,2}(\Omega), \|v\|_{W_{0}^{2,2}(\Omega)}\leq 1$, and write $v=v^{1}+v^{2}$, where $v^{1}\in \mbox{span}\{\omega_{k}\}$ and $(v^{2},\omega_{k})=0\ (k=1,...,m)$. Then from \eqref{E4.4}, we see
$$
<u_{m}'',v>=(u_{m}'',v)_{2}=(u_{m}'',v^{1})_{2}=(f,v^{1})_{2}-\tau(\Delta u_{m},\Delta v^{1})_{2}-\beta(\nabla u_{m},\nabla v^{1})_{2}.
$$
Thus
$$
|<u_{m}'',v>|\leq C(\|f\|_{L^{2}(\Omega)}+\|u_{m}\|_{W^{2,2}(\Omega)}),
$$
here we have used the fact that $\|v^{1}\|_{W^{2,2}(\Omega)}\leq 1$. Consequently
\begin{align}\label{E4.7}
\begin{split}
\int_{0}^{T}\|u_{m}''\|_{W^{-2,2}(\Omega)}dt&\leq C\int_{0}^{T}(\|f\|_{L^{2}(\Omega)}+\|u_{m}\|_{W^{2,2}(\Omega)})dt\\&
\leq C(\|u^{0}\|_{W^{2,2}(\Omega)}+\|u^{1}\|_{L^{2}(\Omega)}+\|f\|_{L^{2}(0,T;L^{2}(\Omega))}).
\end{split}
\end{align}

Now from \eqref{E4.6} and \eqref{E4.7}, we see that there exist a subsequence $\{u_{k}\}_{k=1}^{\infty}$ and $u\in L^{2}(0,T; W_{0}^{2,2}(\Omega))$,
with $u'\in L^{2}(0,T; L^{2}(\Omega)), u''\in L^{2}(0,T; W^{-2,2}(\Omega))$, such that
$$
u_{k}\rightharpoonup u\ \ \ \mbox{in}\ \  L^{\infty}(0,T;W^{2,2}(\Omega))\cap W^{1,\infty}(0,T; L^{2}(\Omega))\cap W^{2,2}(0,T; W^{-2,2}(\Omega)).
$$

Now as in Theorem \ref{T3.1}, we choose a function of the form
$$
h_{j}(x,t)=\sum_{m=1}^{j}\alpha_{j,m}(t)\omega_{m}(x)
$$
such that
$$
\|h_{j}-h\|_{L^{2}(0,T; W_{0}^{2,2}(\Omega))}\to0\  \mbox{as}\ \  j\to\infty
$$
for some $h\in C^{1}([0,T]; C_{0}^{2}(\Omega))$. Here $\{\alpha_{j,m}(t)\}_{m=1}^{k}$ are given smooth functions.  Now multiplying the first equation of (\ref{E4.4}) by
$\alpha_{j,m}$ and summing on $m$ from 1 up to $j$, we, by taking the limit for $k\to\infty$, see that
$$
\int_{0}^{T}<u''(t),h_{j}>+\int_{0}^{T}\beta(\Delta u,\Delta h_{j} )_{2}+\int_{0}^{T}\tau(\nabla u, \nabla h_{j} )_{2}=\int_{0}^{T}(f(t),h_{j})_{2}.
$$
Letting $j\to\infty$ then we are led to
\begin{equation}\label{E4.8}
\int_{0}^{T}<u''(t),h>+\int_{0}^{T}\beta(\Delta u,\Delta h)_{2}+\int_{0}^{T}\tau(\nabla u, \nabla h)_{2}=\int_{0}^{T}(f(t),h)_{2}.
\end{equation}
Since $C^{1}([0,T]; C_{0}^{2}(\Omega))$ is dense in $L^{2}(0,T; W_{0}^{2,2}(\Omega))$,  we conclude equality \eqref{E4.8} is valid for any $h\in L^{2}(0,T; W_{0}^{2,2}(\Omega))$, which further implies
$$
<u''(t),h>+\beta(\Delta u,\Delta h)_{2}+\tau(\nabla u, \nabla h)_{2}=(f(t),h)_{2}
$$
for all $h\in W_{0}^{2,2}(\Omega)$ and a.e. $0\leq t\leq T$. Using the same argument as Theorem \ref{T3.1}, we can also prove $u(0)=u^{0}, u'(0)=u^{1}$, here we omit its details.  Hence
$$
u\in L^{\infty}(0,T;W^{2,2}(\Omega))\cap W^{1,\infty}(0,T; L^{2}(\Omega))\cap W^{2,2}(0,T; W^{-2,2}(\Omega))
$$
is a weak solution of \eqref{E4.1}. Besides, we note that $$u'(t)\in C([0,T]; W^{-2,2}(\Omega))\cap L^{\infty}(0,T; L^{2}(\Omega))$$ and $ L^{2}(\Omega)\subset W^{-2,2}(\Omega)$ with a continuous injection, we, by Lemma \ref{L2.3}, have
$u'(t)\in C_{w}([0,T];L^{2}(\Omega))$. Similar, we also have $u\in C_{w}(0,T;W^{2,2}(\Omega))$. The uniqueness follows from the standard contradiction argument.

\medskip
\begin{thm}\label{T4.2}(Improved regularity).
If $u^{0}\in W^{4,2}(\Omega)\cap W^{2,2}_{0}(\Omega)\ (\mbox{or}\  W^{4,2}(\Omega)\cap W^{1,2}_{0}(\Omega)), u^{1}\in W^{2,2}_{0}(\Omega)\ (\mbox{or}\ W^{2,2}(\Omega)\cap W^{1,2}_{0}(\Omega))
 \ \ f'\in L^{2}(0,T; L^{2}(\Omega)),$ then
\begin{align}\label{E4.9}
\begin{split}
&u\in C_{w}(0,T; W^{4,2}(\Omega)), u'\in C_{w}(0,T; W^{2,2}_{0}(\Omega))\ (\mbox{or}\  C_{w}(0,T; W^{2,2}(\Omega)\cap W^{1,2}_{0}(\Omega))\\
&u''\in L^{\infty}(0,T; L^{2}(\Omega)).
\end{split}
\end{align}
with the estimate
\begin{align}\label{E4.10}
\begin{split}
&\mbox{ess}\sup_{0\leq t\leq T}(\|u'(t)\|_{W^{2,2}(\Omega)}^{2}+\|u''(t)\|_{2}^{2}+\|u(t)\|_{W^{4,2}(\Omega)}^{2})\\&
\leq C(\|f\|_{W^{1,2}(0,T;L^{2}(\Omega))}^{2}+\|u^{0}\|^{2}_{W^{4,2}(\Omega)}+\|u^{1}\|^{2}_{W^{2,2}(\Omega)}).
\end{split}
\end{align}
Here the constant $C$ depends only on $\Omega, T,\beta,\tau$.
\end{thm}

{\bf Proof.} Fix a positive integer $m$ and write $\tilde{u}_{m}:=u_{m}'$, we obtain by differentiating the identity \eqref{E4.4} with respect to $t$,
$$
(\tilde{u}_{m}'',\omega_{k})_{2}+\beta(\Delta\tilde{u}_{m},\Delta\omega_{k})_{2}+\tau(\nabla\tilde{u}_{m},\nabla\omega_{k})_{2}=(f',\omega_{k})_{2}.
$$
Multiplying by $\frac{d^{2}}{dt^{2}}g_{m}^{k}$ and adding for $k=1,...,m$, we discover
$$
(\tilde{u}_{m}'',\tilde{u}_{m}')_{2}+\beta(\Delta\tilde{u}_{m},\Delta\tilde{u}_{m}')_{2}+\tau(\nabla\tilde{u}_{m},\nabla\tilde{u}_{m}')_{2}=(f',\tilde{u}_{m}')_{2}.
$$
and then
\begin{align}\label{E4.11}
\frac{d}{dt}\Big(\|\tilde{u}_{m}'\|_{2}^{2}+\beta\|\Delta\tilde{u}_{m}\|_{2}^{2}+\tau\|\nabla\tilde{u}_{m}\|_{2}^{2}\Big)\leq C(\|\tilde{u}_{m}'\|_{2}^{2}+\|f'\|_{2}^{2}).
\end{align}
Now write
$$
\eta(t):=\|\tilde{u}_{m}'\|_{L^{2}(\Omega)}^{2}+\beta\|\Delta \tilde{u}_{m}\|_{L^{2}(\Omega)}^{2}+\tau\|\nabla \tilde{u}_{m}\|_{L^{2}(\Omega)}^{2}.
$$
Then inequality \eqref{E4.11} reads
$$
\eta'(t)\leq C(\eta(t)+\|f'\|_{L^{2}}^{2})
$$
for $0\leq t\leq T$.

Besides, we note
\begin{align}\label{E4.12}
(f-u''_{m}(t),\omega_{k})_{2}=\beta(\Delta^{2}u_{m},\omega_{k})_{2}-\tau(\Delta u_{m},\omega_{k})_{2}.
\end{align}
Multiplying \eqref{E4.12} by $\lambda_{k}g_{m}^{k}(t)$ and summing $k=1,...,m$, we deduce
\begin{align}\label{E4.13}
\begin{split}
\|\Delta^{2}u_{m}\|_{2}^{2}&\leq(f-u_{m}'',(\beta\Delta^{2}-\tau\Delta)u_{m})_{2}+C\|u_{k}\|_{W^{2,2}(\Omega)}\\&
\leq C(\|f\|_{2}^{2}+\|u''_{m}(t)\|_{2}^{2}+\|u_{m}\|^{2}_{W^{2,2}(\Omega)}).
\end{split}
\end{align}
Applying Grownwall's inequality, we have
\begin{align}\label{E4.14}
\eta(t)\leq e^{Ct}\Big(\eta(0)+\int_{0}^{t}\|f'(s)\|_{L^{2}}^{2}ds\Big),
\end{align}
where
\begin{align*}
\eta(0)&=\|u''_{m}(0)\|_{L^{2}}^{2}+\beta\|\Delta u'_{m}(0)\|_{L^{2}}^{2}+\tau\|\nabla u'_{m}(0)\|_{L^{2}}^{2}\\&
\end{align*}
Employing \eqref{E4.4} and the fact
\begin{align*}
\|u_{k}(0)\|_{W^{4,2}(\Omega)}^{2}\leq C\|u^{0}\|_{W^{4,2}(\Omega)}^{2}.
\end{align*}
we have
\begin{align}\label{E4.15}
\eta(0)\leq C(\|u^{0}\|_{W^{4,2}(\Omega)}^{2}+\|u^{1}\|_{W^{2,2}(\Omega)}^{2}).
\end{align}
Combining \eqref{E4.13}-\eqref{E4.15}, we, by passing to limits as $m=m_{l}\to\infty$, obtain \eqref{E4.10}. Finally, we deduce \eqref{E4.9} by Lemma \ref{L2.2}.

\vskip 0.2in

{\bf Proof of Theorem \ref{T1.2}.}\  As in the proof of Theorem \ref{T1.1}, we only consider the Dirichlet boundary condition. Now define
$$
\mathcal{X}_{T}:=L^{\infty}(0,T;W^{4,2}(\Omega))\cap W^{1,\infty}(0,T; W_{0}^{2,2}(\Omega))\cap W^{2,\infty}(0,T;L^{2}(\Omega))
$$
with norm
$$
\|v\|_{\mathcal{X}_{T}}^{2}:=\mbox{ess}\sup_{0\leq t\leq T}(\|v\|_{W^{4,2}(\Omega)}^{2}+\|v_{t}\|_{W^{2,2}(\Omega)}^{2}+\|v_{tt}\|_{2}^{2}).
$$
And define
\begin{align*}
\bar{M}(R,T):=\{v\in  \mathcal{X}_{T}:  \|v\|_{\mathcal{X}_{T}}\leq R\}
\end{align*}
with $R$ satisfying $C(\Omega)R<1$. Here $C(\Omega)$ is defined in \eqref{E3.15}.  From \eqref{E3.15},  we have
\begin{equation} \label{E4.16}
u(t)\in \bar{M}(R,T)\Rightarrow \|u\|_{L^{\infty}(\Omega\times(0,T))}\leq C(\Omega)R<1,
\end{equation}
and further implies $\frac{1}{(1-u)^{2}}\in W^{1,2}(0,T; L^{2}(\Omega))$.

Now let $r\in (0, R)$ be fixed and
$$
u_{i}(t) \in\bar{M}(r,T),
$$ for $i=1,2$, then by the Theorem \ref{T4.1} and Theorem \ref{T4.2}, the initial-Dirichlet linear problem
\begin{equation}\label{E4.17}
\left \{
\begin{array}{ll}
v_{tt}+ \beta \Delta^{2}v-\tau \Delta v=\frac{\lambda}{(1-u_{i})^{2}},& x\in\Omega, t>0,\\
v(x,0)=u^{0}(x), v_{t}(x,0)=u^{1}(x) & x\in\Omega,   \\
v=\frac{\partial v}{\partial\nu}=0,& x\in \partial\Omega, t>0,\\
 \end{array}
\right.
\end{equation}
has a unique solution $$v_{i}(t):=\mathcal{F}(u_{i})\in L^{\infty}(0,T;W^{4,2}(\Omega))\cap W^{1,\infty}(0,T; W_{0}^{2,2}(\Omega))\cap W^{2,\infty}(0,T;L^{2}(\Omega))$$ for $i=1,2$.

Using the Theorem \ref{T4.2} again, we see that
\begin{align*}%\label{E2.20}
\begin{split}
\|v_{1}-v_{2}\|_{\mathcal{X}_{T}}&\leq\lambda C \|\frac{1}{(1-u_{1})^{2}}-\frac{1}{(1-u_{2})^{2}}\|_{W^{1,2}(0,T;L^{2}(\Omega))}\\
&=\lambda C\Big(\int_{0}^{T}\int_{\Omega}|\frac{1}{(1-u_{1})^{2}}-\frac{1}{(1-u_{2})^{2}}|^{2}dxdt\Big)^{\frac{1}{2}}\\&+2\lambda C\Big(\int_{0}^{T}\int_{\Omega}|\frac{u_{1}'}{(1-u_{1})^{3}}-\frac{u_{2}'}{(1-u_{2})^{3}}|^{2}dxdt\Big)^{\frac{1}{2}}\\
&=:I+II.
\end{split}
\end{align*}
For $I$, we have
\begin{align}\label{E4.18}
\begin{split}
I&\leq 2\lambda C \Big(\int_{0}^{T}\int_{\Omega}\frac{(u_{1}-u_{2})^{2}}{(1-(\theta u_{1}+(1-\theta)u_{2}))^{6}}dxdt\Big)^{\frac{1}{2}}
\\&\leq 2\lambda C k(r)\Big(\int_{0}^{T}\int_{\Omega}(u_{1}-u_{2})^{2}dxdt\Big)^{\frac{1}{2}}\\&
\leq 2\lambda T^{\frac{1}{2}} C k(r)\|u_{1}-u_{2}\|_{\mathcal{X}_{T}}.
\end{split}
\end{align}
For $II$, we have
\begin{align}\label{E4.19}
\begin{split}
II&\leq 2\lambda C \Big(\int_{0}^{T}\int_{\Omega}\frac{(u'_{1}-u'_{2})^{2}}{(1-u_{1})^{6}}dxdt\Big)^{\frac{1}{2}}\\&
+2\lambda C\Big(\int_{0}^{T}\int_{\Omega}|u'_{2}|^{2}\Big|\frac{1}{(1-u_{1})^{3}}-\frac{1}{(1-u_{2})^{3}}\Big|^{2}dxdt\Big)^{\frac{1}{2}}\\&
\leq2\lambda CT^{\frac{1}{2}}k(r)\Big(\|u_{1}-u_{2}\|_{W^{1,\infty}(0,T;L^{2}(\Omega))}+\|u_{1}-u_{2}\|_{L^{\infty}((0,t)\times\Omega)}\|u_{2}\|_{W^{1,\infty}(0,T;L^{2}(\Omega))}\Big)\\&
\leq2\lambda C(r+k(r))T^{\frac{1}{2}}\|u_{1}-u_{2}\|_{\mathcal{X}_{T}}.
\end{split}
\end{align}
Here and in what follows $k(r)$ is a positive nondecreasing function for $r\in [0,R_{0}]$ and $C$  depends only on $\Omega, T,\beta,\tau$. From \eqref{E4.18} and \eqref{E4.19}, we have
\begin{align}\label{E4.20}
\|v_{1}-v_{2}\|_{\mathcal{X}_{T}}\leq 2\lambda T^{\frac{1}{2}} C(r+k(r))\|u_{1}-u_{2}\|_{\mathcal{X}_{T}}.
\end{align}

Now consider the unique solution $w(t)$ to the linear problem
$$
w_{tt}+\beta \Delta^{2}w-\tau \Delta w=0,\ \ \ x\in \Omega, t>0,
$$
with the same boundary and initial conditions as \eqref{E4.1}.  Obviously we have, by the Theorem \ref{T4.2}
 $$w(t)\in L^{\infty}(0,T;W^{4,2}(\Omega))\cap W^{1,\infty}(0,T; W_{0}^{2,2}(\Omega))\cap W^{2,\infty}(0,T;L^{2}(\Omega))$$
  such that
 \begin{align}\label{E4.21}
 \begin{split}
\|w\|_{\mathcal{X}_{T}}&\leq C(\|u^{0}\|_{W^{4,2}(\Omega)}+\|u^{1}\|_{W^{2,2}(\Omega)})=C\rho
\end{split}
\end{align}

Define the ball
$$
\mathcal{B}_{\frac{r}{2}}=\{u\in \mathcal{X}_{T}: \|u-w\|_{\mathcal{X}_{T}}\leq \frac{r}{2}\}.
$$
Choosing $\rho$  small enough  such that
$$
C\rho+\frac{r}{2}<r,
$$
we then have, if $u(t)\in \mathcal{B}_{\frac{r}{2}}$,
$$
\|u(t)\|_{\mathcal{X}_{T}}\leq \|w(t)\|_{\mathcal{X}_{T}}+\frac{r}{2}<r,\ \ \forall\ 0\leq t\leq T,
$$
and then $\mathcal{B}_{\frac{r}{2}}\subset M(r,T)$. \medskip

Now using estimate \eqref{E4.20}, we find
\begin{align}\label{E4.22}
\|v_{i}-w\|_{\mathcal{X}_{T}}\leq 2\lambda T^{\frac{1}{2}} C(r+k(r))r
\end{align}
for $i=1,2$. Now choosing $\lambda$ so small that
$$
\lambda\leq \lambda(r,T):= \frac{1}{4T^{\frac{1}{2}}C(k(r)+r)r},
$$
we have from \eqref{E4.20} and \eqref{E4.22} that
\begin{align*}
\begin{split}
&\|v_{1}-v_{2}\|_{\mathcal{X}_{T}}\leq \frac{1}{2}\|u_{1}-u_{2}\|_{\mathcal{X}_{T}};\\
&\|v_{i}-w\|_{\mathcal{X}_{T}}\leq\frac{r}{2}.
\end{split}
\end{align*}
Hence the map
\begin{align}
\begin{split}
\mathcal{F}:& \mathcal{B}_{\frac{r}{2}}\to \mathcal{B}_{\frac{r}{2}}\\
& u_{i}\to v_{i}\ \ (i=1,2)
\end{split}
\end{align}
is a contraction map and it has a unique fixed point $u=\mathcal{F}(u)$ in $\mathcal{B}_{\frac{r}{2}}$ for $0<\lambda\leq\lambda(r,T)$. Finally, we claim that this solution satisfies
$$
u(t)\in C_{w}([0,T]; W^{4,2}(\Omega)), u'(t)\in  C_{w}([0,T]; W_{0}^{2,2}(\Omega)).
$$
Indeed, since $u'(t)\in C([0,T]; L^{2}(\Omega))\cap L^{\infty}(0,T; W_{0}^{2,2}(\Omega))$ and $W_{0}^{2,2}(\Omega)\subset L^{2}(\Omega)$ with a continuous injection, and then $u'(t)\in  C_{w}([0,T]; W_{0}^{2,2}(\Omega))$ by Lemma \ref{L2.2}. Similar, we have $u(t)\in C_{w}([0,T]; W^{4,2}(\Omega))$.

The proof of (ii) is similar with (iii) of Theorem \ref{T1.1}, we omit it here.

\end{document}